\documentclass[12pt]{amsart}

\RequirePackage{hyperref}

\theoremstyle{plain}
\newtheorem{prop}{Proposition}

\newcommand\Aut{\operatorname{Aut}}
\newcommand\pr{\noindent\textit{Proof} : }

\newcommand\Z{\mathbb{Z}}

\newcommand\C{\mathbb{C}}
\newcommand\Tr{\operatorname{Tr}}
\def\lr#1{\langle {#1} \rangle}

\begin{document}

\title[Quotients of Jacobians]{Quotients of Jacobians\protect}
\author[Arnaud Beauville]{Arnaud Beauville}
\address{Universit\'e C\^ote d'Azur\\
CNRS -- Laboratoire J.-A. Dieudonn\'e\\
Parc Valrose\\
F-06108 Nice cedex 2, France}
\email{arnaud.beauville@unice.fr}

\begin{abstract}
Let $C$ be a  curve of genus $g$, and let $G$ be a finite group of automorphisms of $C$. The group $G$ acts on the Jacobian $J$ of $C$; we prove that for $g\geq 21$ the quotient $J/G$  has canonical singularities and Kodaira dimension $0$.  On the other hand we give examples with $g\leq 4$ for which   $J/G$ is uniruled. 
\end{abstract}

\maketitle


\begin{center}
\emph{Pour Enrico, apr\`es 50 ans d'amiti\'e}
\end{center}
\medskip	
\section{Introduction}
Let $C$ be a (smooth, projective, complex) curve of genus $g$, and let $G$ be a finite group of automorphisms of $C$. The group $G$ acts on the Jacobian $J$ of $C$. There are two possibilities for the quotient $J/G$ \cite{K-L}\footnote{I am indebted to A. H\"oring for pointing out this reference}: either it has canonical singularities and Kodaira dimension $0$, or it is uniruled. The aim of this note is to show that the latter case is rather exceptional -- in fact, it does not occur for $g\geq 21$. This bound is  rough, and can certainly be lowered. On the other hand we give examples of uniruled $J/G$ for $g\leq 4$. 
When $g=3$ this is an essential ingredient in our proof that the cycle $[C]-[(-1)^*C]$ on $J$ is torsion modulo algebraic equivalence \cite{B-S}; in fact this note grew out of the observation that already in genus 3, the case where $J/G$ is uniruled is quite rare. 

\section{$J/G$ canonical}

Let $C$ be a curve of genus $g$ and $J$ its Jacobian. Let $G$ be a subgroup of $\Aut(C)$, hence also of $\Aut(J)$. By \cite[Theorem 2]{K-L}, there are two possibilities:

$\bullet$ Either $J/G$ has canonical singularities, and Kodaira dimension $0$;

$\bullet$ or $J/G$ is uniruled.

\begin{prop}
Assume $g\geq 21$. The quotient variety $J/G$ has canonical singularities, and Kodaira dimension $0$.
\end{prop} 
\pr 
We first observe that the fixed locus  $\operatorname{Fix}(\sigma )\subset  J$ of any element $\sigma \neq 1$ of $G$ 
has codimension $\geq 2$. Indeed the dimension of  $\operatorname{Fix}(\sigma )$  is the multiplicity of the eigenvalue $1$ for the action of $\sigma $ on $H^0(C,K_C)$, that is, the genus of $C/\lr{\sigma }$, which is $\leq g-2$   as soon as $g\geq 4$ by the Hurwitz formula.

Let $\sigma $ be an element of order $r$ in $G$, and let $p\in J$ be a fixed point of $\sigma $. The action of $\sigma $ on the tangent space $T_p(J)$ is isomorphic to the action on $T_0(J)=H^{0,1}(C)$.   
Let $\zeta =e^{\frac{2\pi i}{r} }$; we write the eigenvalues of $\sigma $ acting on $H^{0,1}(C)$ in the form $\zeta ^{a_1},\ldots , \zeta ^{a_g}$, with $0\leq a_i<r$. By Reid's criterion \cite[Theorem 3.1]{R}, $J/G$ has canonical singularities if and only if $\sum a_i\geq r$ for all $\sigma $ in $G$. 

The eigenvalues of $\sigma $ acting on $H^1(C,\C)$ are $\zeta ^{a_1},\ldots , \zeta ^{a_g}; {\zeta }^{-a_1},\ldots , {\zeta} ^{-a_g}$.  
Thus $\Tr \sigma ^*_{|H^1(C,\C)}=\allowbreak 2\sum \cos  \dfrac{2\pi a_i}{r} \,\cdot  $ By the Lefschetz formula, this trace is equal to $2-f$, where $f$ is the number of fixed points of $\sigma $; in particular, it is $\leq 2$. Using $\cos x\geq 1-\dfrac{x^2}{2} $, we find
\[1\geq  \sum \cos \dfrac{2\pi a_i}{r}\geq  g-\dfrac{ 2\pi ^2}{r^2} \sum a_i^2\ ,\]
hence $\displaystyle\ (\sum a_i)^2\geq \sum a_i^2 \geq \dfrac{g-1}{2 \pi ^2}\, r^2$. Thus if 
$g\geq 1+2\pi ^2=20.739\ldots $, we get $\sum a_i>r$.\qed

\medskip	
\noindent\emph{Remark}$.-$ The proposition does not extend to the case of an arbitrary abelian variety: indeed we give below examples where $J/G$ is uniruled; then if $A$ is any abelian variety, the quotient of $J\times A$ by $G$ acting trivially on $A$
is again uniruled.

\section{$J/G$ uniruled}

We will now give examples of (low genus) curves $C$ such that $J/G$ is uniruled.
The genus 2 case is quite particular (and probably well known). We put $\rho :=e^{\frac{2\pi i}{3} }$. We assume $G\neq \{1\} $.
\begin{prop}
If $g(C)=2$, $J/G$ has Kodaira dimension $0$ if and only if $G=\lr{\sigma }$ and $\sigma $ is the hyperelliptic involution, or an automorphism of order $3$ with eigenvalues $(\rho ,\rho ^2)$ on $H^0(C,K_C)$, or
an automorphism of order $6$ with eigenvalues $(-\rho ,-\rho ^2)$ on $H^0(C,K_C)$.

\end{prop}
\pr  Let $\sigma $ be an element of order $r$ in $G$, and let $\zeta $ be a primitive $r$-th root of unity. 
The eigenvalues of $\sigma $ acting on $H^0(C,K_C)$ are $\zeta^{a} ,\zeta ^{b}$, with  $0\leq a,b<r$. Suppose that $J/G$ has canonical singularities.  By Reid's criterion we must have $a+b\geq r$; replacing $\zeta $ by $\zeta ^{-1}$ gives $2r-a-b\geq r$, hence $a+b=r$. Since $\sigma $ has order $r$, $\sigma ^k$ acts non-trivially on $H^0(C,K_C)$ for $0<k<r$; therefore $a$ and $r$ are coprime. Replacing $\zeta $ by $\zeta ^{a}$ we may assume that the eigenvalues of $\sigma $ are $\zeta $ and $\zeta ^{-1}$.

Put $\zeta =e^{i\alpha }$; the trace of $\sigma $ acting on $H^1(C,\Z)$ is $2(\zeta +\zeta ^{-1})=4 \cos \alpha $, and it is an integer. This implies $\cos \alpha =0$ or $\pm \dfrac{1}{2}$, hence $\alpha \in \Z\cdot \dfrac{2\pi }{6}  $ and $r\in\{2,3,6\} $. If $r=2$ we must have $\zeta =\zeta ^{-1}=-1$, hence $\sigma $ is the hyperelliptic involution; if $r=3$ or $6$ we have $\{\zeta ,\zeta ^{-1}\}=\{\rho ,\rho ^2\}  $ or $\{-\rho ,-\rho ^2\} $.

Thus $G$ contains elements of order $3$ and $6$, and at most one element of order $2$.
Since the order of $\Aut(C)$  is not divisible by $9$
 \cite[11.7]{B-L},  $G$ is isomorphic to either $\Z/3$, or a central extension of $\Z/3$ by $\Z/2$ -- that is, $\Z/6$. This proves the Proposition.
\qed

\medskip	
The last two  cases are realized by the curves $C$ of the form $y^2=x^6+kx^3+1$, with the automorphisms $\sigma _1:(x,y)\mapsto (\rho x, y)$ and $\sigma _2:(x,y)\mapsto (\rho x,-y)$. The surfaces $S_i:=J/\lr{\sigma _i}$ are (singular) K3 surfaces, each depending on one parameter. The surface $S_1$ is studied in detail in \cite{B-V}: it has 9 singular points of type $A_2$. The surface $S_2$ is the quotient of $S_1$ by the involution induced by $(-1)_J$; it has one singular point of type $A_5$ (the image of $0\in J$), 4 points of type $A_2$ corresponding to pairs of singular points of $S_1$, and $5$ points of type $A_1$ corresponding to triples of points of order 2 of $J$.
\bigskip	

It is more subtle to give examples with $g\geq 3$. By Reid's criterion, we must exhibit an automorphism $\sigma $ of  $C$, of order $r$, such that the eigenvalues of $\sigma $ acting on $H^0(C,K_C)$ are of the form $\zeta ^{a_1},\ldots ,\zeta ^{a_g}$, with $\zeta $ a primitive $r$-th root of unity, $a_i\geq 0$ and $\sum a_i<r$. In the following tables we give the curve $C$ (in affine coordinates), the order $r$ of $\sigma $ and $\zeta $, the expression of $\sigma $, a basis of eigenvectors for $\sigma $ acting on $H^0(C,K_C)$, the exponents $a_1,\ldots ,a_g$; in the last column we check the inequality $\sum a_i<g$, which implies that $J/\lr{\sigma }$ is uniruled.

\vfill\eject
\begin{center}
\framebox{$g=3$}
\end{center}
{\renewcommand{\arraystretch}{2}
\[\begin{array}{|c|c|c|c|c|c|}
        \hline
         \text{curve} &  r & \sigma(x,y)= &   \text{basis of } H^0(K_C)  & a_1,\ldots ,a_3  &  \sum a_i\\ 
          \hline\rule[-0.5cm]{0cm}{2pt}
         y^2=x^8-1&  8 & (\zeta x,  y)& \dfrac{dx}{y}, \dfrac{xdx}{y}, \dfrac{x^2dx}{y}& 1,2,3 &  6<8\\
          \hline\rule[-0.5cm]{0cm}{2pt}
       y^2=x(x^6-1) &  12 & (\zeta ^2x,\zeta y)& \dfrac{dx}{y}, \dfrac{xdx}{y}, \dfrac{x^2dx}{y}& 1,3,5 &  9<12\\
        \hline\rule[-0.5cm]{0cm}{2pt}
       y^2=x(x^7-1) &  14 & (\zeta ^2x,\zeta y)& \dfrac{dx}{y}, \dfrac{xdx}{y}, \dfrac{x^2dx}{y}& 1,3,5 &  9<14\\   
        \hline\rule[-0.5cm]{0cm}{2pt}
       y^3=x(x^3-1) &  9 & (\zeta ^3x,\zeta y)& \dfrac{dx}{y^2}, \dfrac{xdx}{y^2}, \dfrac{dx}{y}& 1,4,2 &  7<9\\   
   \hline\rule[-0.5cm]{0cm}{2pt}
       y^3=x^4-1 &  12 & (\zeta ^{-3}x,\zeta^{-4} y)& \dfrac{dx}{y^2}, \dfrac{xdx}{y^2}, \dfrac{dx}{y}& 5,2, 1 & 8<12\\   

\hline
\end{array} \]
\medskip	

The fact that $J/\lr{\sigma }$ is uniruled for the fourth curve $C$ is an important ingredient of the proof in \cite{B-S} that the Ceresa cycle $[C]-[(-1)^*C]$ in $J$ is torsion modulo algebraic equivalence. The same property for the last curve is observed in \cite{L-S}.  The second curve has been provided to me by D. Conti, A. Ghigi and R. Pignatelli, using the database \cite{CGP}}.
\bigskip	
\begin{center}
\framebox{$g=4$}
\end{center}
{\renewcommand{\arraystretch}{2}
\[\begin{array}{|c|c|c|c|c|c|}
        \hline
         \text{curve} &  r & \sigma(x,y)= &   \text{basis of } H^0(K_C)  & a_1,\ldots ,a_4  & \sum a_i\\ 
          \hline\rule[-0.5cm]{0cm}{2pt}
         y^2=x(x^9-1)&  18 & (\zeta^2 x, \zeta  y)& \dfrac{dx}{y}, \dfrac{xdx}{y}, \dfrac{x^2dx}{y}, \dfrac{x^3dx}{y}& 1,3,5,7 &  16<18\\
        \hline\rule[-0.5cm]{0cm}{2pt}
       y^3=x^5-1 &  15 & (\zeta ^{-3}x,\zeta ^{-5} y)& \dfrac{dx}{y^2}, \dfrac{xdx}{y^2}, \dfrac{x^2dx}{y^2}, \dfrac{dx}{y}& 7,4,1,2 &  14<15\\   
\hline

\end{array} \]

\bigskip	
In genus 5 there is no example:
\begin{prop}
Assume $g=5$. The quotient variety $J/G$ has canonical singularities, and Kodaira dimension $0$.
\end{prop} 
\pr The paper \cite{K-K} contains the list of possible automorphisms of genus 5 curves, and gives for each automorphism of order $r$  its eigenvalues $\zeta ^{a_1},\ldots ,\zeta ^{a_5} $ on $H^0(K)$, with $\zeta =e^{\frac{2\pi i }{r}}$ and $0\leq a_i<5$. We want to prove $\sum a_i\geq r$.

We first observe that 
if an eigenvalue and its inverse appear in the list, we have $\sum a_i\geq r$; this eliminates all cases with $r\leq 10$. For $r\geq 11$ the numbers
$a_1,\ldots ,a_5$ are always distinct and nonzero; therefore $\sum a_i\geq 1+2+\cdots +5\allowbreak = 15$. The remaining cases are $r=20$ and $22$; one checks immediately that for all $k$ coprime to $r$ we have $\sum b_i\geq r$, where $b_i\equiv ka_i\pmod{r}$ and $0<b_i<r$.  This proves the Proposition.\qed

\medskip	
\noindent \textbf{Added in proof :} R. Serova just proved that there are in fact no examples
in genus  $\geq 6$ (\textsl{Quotients of Jacobians}, preprint \texttt{arXiv:2407.00860}).

\end{document}